\documentclass[12pt]{article}
\usepackage{amsfonts, amsmath, amssymb, latexsym, url}
\usepackage{graphicx}

\newcommand{\includefig}[3]{\includegraphics[scale=#3]{#1.pdf}}

\newcommand{\R}{\ensuremath{\mathbb{R}}}
\newcommand{\qed}{\hfill $\Box$ }

\newcommand{\opstyle}[1]{\ensuremath{\operatorname{#1}}}
\newcommand{\conestyle}[1]{\mbox{\sl #1}}
\newcommand{\rank}{\opstyle{rank}}
\newcommand{\Inc}{\opstyle{Inc}}
\newcommand{\Sk}{\opstyle{Sk}}
\newcommand{\Ri}{\opstyle{Ri}}
\newcommand{\Aut}{\opstyle{Aut}}
\newcommand{\Sym}{\opstyle{Sym}}
\newcommand{\Hyp}{\opstyle{Hyp}}
\newcommand{\OHyp}{\opstyle{OHyp}}
\newcommand{\ONeg}{\opstyle{ONeg}}
\newcommand{\sWPMET}{\conestyle{sWPMET}}
\newcommand{\sPMET}{\conestyle{sPMET}}
\newcommand{\wPMET}{\conestyle{wPMET}}
\newcommand{\PMET}{\conestyle{PMET}}

\newcommand{\sWQMET}{\conestyle{sWQMET}}
\newcommand{\wWQMET}{\conestyle{wWQMET}}
\newcommand{\WQMET}{\conestyle{WQMET}}
\newcommand{\QMET}{\conestyle{QMET}}
\newcommand{\dWMET}{\conestyle{dWMET}}
\newcommand{\sWMET}{\conestyle{sWMET}}
\newcommand{\WMET}{\conestyle{WMET}}
\newcommand{\MET}{\conestyle{MET}}

\newcommand{\OMCUT}{\conestyle{OMCUT}}
\newcommand{\OCUT}{\conestyle{OCUT}}
\newcommand{\WCUT}{\conestyle{WCUT}}
\newcommand{\CUT}{\conestyle{CUT}}

\newcommand{\PHYP}{\conestyle{PHYP}}
\newcommand{\QHYP}{\conestyle{QHYP}}
\newcommand{\WHYP}{\conestyle{WHYP}}
\newcommand{\HYP}{\conestyle{HYP}}
\newcommand{\OWHYP}{\conestyle{OWHYP}}

\newtheorem{theor}{Theorem}

\newtheorem{lemma}{Lemma}

\newtheorem{conj}{Conjecture}

\begin{document}

\title{{\bf Cones of Weighted and Partial Metrics}}

\author{Michel Deza\footnote{michel.deza@ens.fr, \'Ecole Normale
Sup\'erieure, Paris}~ 
 Elena Deza\footnote{elena.deza@gmail.com, Moscow State Pedagogical
University, Moscow}~ and
  Jano\v{s} Vidali\footnote{janos.vidali@fri.uni-lj.si, University of 
Ljubljana, Slovenia}~
}
\date{}

\maketitle

\begin{abstract}
\noindent A {\em partial semimetric} on  $V_{n}=\{1, \dots, n\}$ 
is  a function $f=((f_{ij})):   
V_n^2\longrightarrow\mathbb{R}_{\ge 0}$
satisfying 
 $f_{ij}=f_{ji}\ge f_{ii}$ and 
$f_{ij}+f_{ik}-f_{jk}-f_{ii} \ge 0$ 
for all $i,j,k\in V_n$.
The function $f$ is a {\em weak partial semimetric} if 
 $f_{ij}\ge f_{ii}$ is dropped,  
 and it is a {\em strong 
partial semimetric} 
if
 $f_{ij}\ge f_{ii}$ is complemented by   $f_{ij}\le f_{ii}+f_{jj}$.

We describe the cones of weak and strong partial semimetrics 
via corresponding weighted semimetrics and list their $0,1$-valued 
elements, identifying when they belong to extreme rays. We consider also 
related cones, including those of {\em partial hypermetrics}, {\em 
weighted hypermetrics}, {\em $\ell_1$-quasi semimetrics} and 
weighted/partial cuts.

\end{abstract}


\section{Convex cones under consideration}

There are following two  main motivations for this study.
One is to extend the rich theory  of 
metric, cut and  hypermetric cones on weighted, partial and non-symmetric
generalizations of metrics. 
Another is a  new  appoach to  partial semimetrics (having important applications in Computer 
Science) via cones formed by them.

A {\em convex cone in $\R^m$} (see, for example,~\cite{Sch86}) is defined
either
by {\em generators} $v_1$, \dots, $v_{N}$, as $\{ \sum \lambda_i v_i:
\lambda_i\geq 0\}$, or by {\em linear inequalities} $f^{1}$, \dots, 
$f^{M}$,
as $\cap_{j=1}^M \{x\in \R^m: f^j(x)=\sum_{i=1}^mf_i^jx_i\geq 0\}$. 

Let $C$ be an $m'$-dimensional convex cone in $\R^{m}$. 
Given $f\in \R^{m}$, the linear inequality 
$f(x)=\sum_{i=1}^mf_ix_i=\langle f,x \rangle  \geq 0$ is 
said to
be {\em valid} for $C$ if it holds for all $x \in C$. Then the set
$\{ x \in C: \langle f,x \rangle= 0 \}$ is called the {\em face}
of $C$, {\em induced by} $F$.
A face of dimension $m'-1$, $m'-2$, $1$ is  
called a {\em
facet}, {\em ridge}, {\em extreme ray}  of $C$, respectively
(a {\em ray} is a set $\R_{\ge 0} x$ with $x\in C$).
Denote by $F(C)$ the set of facets of $C$ and by $R(C)$ the set of its 
extreme rays.
We consider only {\em polyhedral} (i.e., $R(C)$ and, alternatively, $F(C)$ 
is finite) 
{\em pointed} (i.e., $(0)\in C$) convex 
cones.
Each ray $r\subset C$ below contains a unique  {\em good 
representative}, i.e.,  an integer-valued 
vector $v(r)$ with g.c.d. $1$ of its entries; so, by abuse of language, we 
will 
identify $r$ with $v(r)$.

For a ray $r\subset C$ denote
by $F(r)$ 
 the set $\{f\in F(C): r\subset f\}$. For a face $f\subset C$ denote
by $R(f)$
 the set $\{r\in R(C): r\subset f\}$.   
The {\em incidence number}  $\Inc(f)$ of a face $f$ (or  $\Inc(r)$ of a 
ray 
$r$) is the number
  $|\{r\in R(C): r\subset f\}|$ (or, respectively,
 $|\{f\in F(C): r\subset f\}|$). The $\rank(f)$ of 
 a face
 $f$ (or $\rank(r)$ of 
a ray $r$) is the dimension of 
  $\{r\in R(C): r\subset f\}$ (or 
 of $\{f\in F(C): r\subset f\}$). 

Two extreme rays (or facets) of $C$  are  {\em adjacent on
$C$} if
they span a $2$-dimensional face (or, respectively, their intersection
has dimension $m'-2$).
The {\em skeleton} $\Sk(C)$ is the graph
whose vertices are the extreme rays of $C$ and with an
 edge between two vertices if the corresponding rays are adjacent on
$C$.
The {\em ridge graph} $\Ri(C)$ is the graph 
 whose vertices are facets of $C$ and with an edge between two vertices if
 the corresponding facets are adjacent on
$C$. Let $D(G)$ denote the diameter of the graph $G$

Given a cone  $C_n$ of   some  functions, say,  $d=((d_{ij})):
V_n^2\longrightarrow\mathbb{R}_{\ge 0}$  
the {\em $0$-extension} of the inequality $\sum_{1\leq i\not=
j\leq n-1} F_{ij}d_{ij}\geq 0$ is the inequality
\begin{equation*}
\sum_{1\leq i\not= j\leq n}F'_{ij}d_{ij}\geq 0\mbox{~with~} F'_{ni}=F'_{in}=
0\mbox{~and~} F'_{ij}=F_{ij}\mbox{,~otherwise}.
\end{equation*}
Clearly, the $0$-extension of any facet-defining inequality of a cone $C_n$
is a valid inequality (usually,  facet-defining) of $C_{n+1}$.
The $0$-extension of an extreme ray is defined similarly.
For any cone $C$ denote by $0,1$-$C$ the cone generated by all extreme
rays of $C$ containing a non-zero $0,1$-valued point.

The  cones $C$ considered here  will be symmetric under permutations and usually
$\Aut(C)=\Sym(n)$. All orbits below are under $\Sym(n)$.
\clearpage

Set $V_{n}=\{1, \dots, n\}$.  The function $f=((f_{ij})):
V_n^2\longrightarrow\mathbb{R}
$ is called {\em weak partial semimetric} 
if the following holds:

(1) $f_{ij}=f_{ji}$ ({\em symmetry}) for  all $i,j \in V_{n}$,

(2) $L_{ij}: f_{ij}\ge 0$ ({\em non-negativity}) for  all $i,j \in V_{n}$,
and

(3) $Tr_{ij,k}: f_{ik}+f_{kj}-f_{ij}-f_{kk}\ge 0$
 ({\em triangle inequality})  for all $i,j,k \in V_{n}$.

Weak partial
semimetrics were introduced in~\cite{He99} as a generalization of 
{\em partial 
semimetrics}  introduced in~\cite{Ma92}.
Clearly, all $Tr_{ij,i}=0$ and   
$Tr_{ii,k}=2f_{ik}-f_{ii}-f_{kk}=Tr_{ij,k}+Tr_{kj,i}$. So, it is sufficient to 
require (2) only for $i=j$ and (3) only for different $i,j,k$.
The  weak
partial semimetrics on $V_{n}$ form a    ${n+1 \choose 
2}$-dimensional convex cone  
with $n$ facets
$L_{ii}$ and $3{n \choose 3}$ facets $Tr_{ij,k}$. Denote this cone  
 by  $\wPMET_n$.
\vspace{0.7mm}

A weak partial semimetric $f$ is called {\em partial semimetric} if it
holds that 

(4) $M_{ij}: f_{ij}-f_{ii}\ge 0$ ({\em small self-distances})
for  all different $i,j \in V_{n}$.

The  
partial semimetrics on $V_{n}$ form a    ${n+1 \choose
2}$-dimensional  subcone, denote it by $\PMET_n$, of $\wPMET_n$. This cone
  has 
$n$ facets
$L_{ii}$, $2{n \choose 2}$ facets $M_{ij,i}$ and $3{n \choose 3}$ facets 
$Tr_{ij,k}$.
Partial metrics were introduced  by Matthews in~\cite{Ma92}
for treatment of partially defined objects in Computer Science; see also~%
\cite{Ma08,Hi01,Se97}. The cone  $\PMET_n$ was considered in~\cite{DeDe}.
\vspace{0.7mm}

A partial semimetric $f$ is called {\em strong partial semimetric} if it
holds that

(5) $N_{ij}: f_{ii}+f_{jj}-f_{ij}\ge 0$ ({\em large 
self-distances})
for  all $i,j \in V_{n}$.

So, $f_{ii}=N_{ij}+M_{ji}\ge 0$, i.e., (5) and (4) imply $L_{ii}$ for all.   
$i$. The strong partial semimetrics on $V_{n}$ form a    ${n+1 \choose
2}$-dimensional  subcone, denote it by $\sPMET_n$, of $\PMET_n$. This 
cone
  has $3{n+1 \choose 3}$ facets: $2{n \choose 2}$ facets $M_{ij}$, ${n 
\choose 2}$ facets 
$N_{ij}$ and $3{n \choose 3}$ facets
$Tr_{ij,k}$.
\vspace{0.7mm}

A partial semimetric $f$ is called {\em semimetric} if it
holds that

(6) $f_{ii}=0$ ({\em reflexivity})
for  all $i \in V_{n}$.

The semimetrics on $V_{n}$ form a   ${n \choose
2}$-dimensional convex cone, denoted  by $\MET_n$, which 
  has $3{n \choose 3}$ facets
$Tr_{ij,k}$ (clearly, $f_{ij}=\frac{Tr_{ij,k}+Tr_{jk,i}}{2}\ge 0$).
This cone is well-known; see, for example,~\cite{DL}  and references 
there. 
\vspace{0.7mm}

The  function $f$ is  {\em quasi-semimetric} if only  
(2), (3), (6)   are required. The  quasi-semimetrics on $V_n$ 
form a   $n(n-1)$-dimensional convex cone, denoted  by $\QMET_n$, which
  has 
$2{n \choose 2}$ facets
$L_{ij}$ and 
$6{n \choose 3}$ facets
$OTr_{ij,k}: f_{ik}+f_{kj}-f_{ij}\ge 0$ ({\em oriented triangle 
inequality}). But other oriented versions of $Tr_{ij,k}$ (for example,  
$f_{ik}+f_{kj}-f_{ji}$) are not valid on $\QMET_n$.  

A quasi-semimetric $f$ is  {\em weightable} if there exist a
(weight) function $w=(w_i): V_n\longrightarrow\mathbb{R}_{\ge 0}$
such that $f_{ij}+w_i=f_{ji}+w_j$
for  all  $i, j \in V_{n}$. 
  Such 
quasi-semimetrics $f$ (or, equivalently, pairs $(f,w)$) on $V_{n}$ form a 
${n+1 \choose 2}$-dimensional cone, 
 denote it by 
 $\WQMET_n$,  with  $2{n \choose 2}$ facets 
$L_{ij}$  and $3 {n \choose 3}$ facets $OTr_{ij,k}$ since, for a 
quasi-semimetric, 
$OTr_{ij,k}=OTr_{ji,k}$ if it is weightable.

A weightable quasi-semimetric $(f,w)$ with all $f_{ij}\le w_j$
 is  
a {\em weightable strong quasi-semimetric}. But if, on the contrary,
(2) is  weakened to 
$f_{ij}+f_{ji}\ge 0$ (so, $f_{ij}<0$ is allowed),
$(f,w)$ is
a {\em weightable weak quasi-semimetric}. Denote by $\sWQMET_n$ and 
$\wWQMET_n$ the corresponding cones.
\vspace{0.7mm}

Let us denote the function  $f$ by $p$, $d$, or $q$ if it is a weak 
partial 
semimetric, semimetric, or  weightable weak  
quasi-semimetric,
respectively.
\vspace{1.3mm}

A {\em weighted semimetric} $(d;w)$ is a semimetric  
$d$ with a weight function $w: V_n \rightarrow \mathbb{R}_{\ge 0}$
on its points.  Denote by $(d;w)$ the matrix $((d'_{ij}))$, $0\le i,j\le 
n$, 
with $d'_{00}=0$, $d'_{0i}=d'_{i0}=w_i$ for $i\in V_n$ and 
$d'_{ij}=d_{ij}$ for $i,j\in V_n$.
The
weighted semimetrics $(d;w)$ on $V_{n}$ form a   ${n+1 \choose
2}$-dimensional convex cone
 with $n$ facets
$w_{i}\ge 0$ and $3{n \choose 3}$ facets $Tr_{ij,k}$. Denote this cone  by
$\WMET_n$. So,  $\MET_n\simeq \{(d;(k,\dots, 
k)):d\in \WMET_n\}$. Also, $\MET_n=\QMET_n\cap \PMET_n$.

Call a  weighted semimetric $(d;w)$ {\em down-} or 
{\em up-weighted} if

(4') $d_{ij}\ge w_i-w_j$, or 

(5') $d_{ij}\le w_i+w_j$ \\
holds (for all distinct $i,j\in V_n$), respectively.
Denote by $\dWMET_n$ 
 the cone of down-weighted 
semimetrics on $V_n$  and by $\sWMET_n$ the cone 
of {\em strongly}, i.e., both, down- and  up-, weighted
semimetrics. So, 
$\sWMET_n= \MET_{n+1}$.

\section{Maps $P,Q$  and semimetrics}

Given a weighted semimetric  $(d;w)$, define the  map $P$ by 
the function $p=P(d;w)$  with
$p_{ij}=\frac{d_{ij}+w_i+w_j}{2}$. Clearly, $P$ is an {\em automorphism}
(invertible linear operator) of  the vector space 
 $\mathbb{R}^{n+1 \choose
2}$, and $(d;w)=P^{-1}(p)$, where the inverse  map 
$P^{-1}$ is  defined by 
 $d_{ij}=2p_{ij}-p_{ii}-p_{jj}$, 
$w_i=p_{ii}$.

Define the map $Q$ by
the function $(q,w)=Q(d;w)$  with
 $q_{ij}
=\frac{d_{ij}-w_i+w_j}{2}$. So, $Q(d;w)=P(d;w)-((1))w$ 
(i.e., $q_{ij}=p_{ij}-p_{ii}$) and  
 $d_{ij}=q_{ij}+q_{ji}$, is  the {\em symmetrization
semimetric} of  $q$.

{\em Example.} Below are given: the semimetric 
$d=2\delta(\{56\},\{1\},\{23\},\{4\})$ $-\delta(\{56\})\in 
\MET_6$, and, taking weight $w=(1_{i\in \{56\}})=(0,0,0,0,1,1)$, the 
partial semimetric 
$P(d;w)=J(\{56\})+\delta(\{56\},\{1\},\{23\},\{4\})$ (its ray is     
extreme in $\PMET_6$)  
and the weightable 
quasi-semimetric 
$Q(d;w)=\delta'(\{1\})+\delta'(\{23\})+\delta'(\{4\})$ (its ray is not 
extreme in $\WQMET_6$).

\clearpage

\begin{center}

${\bf 0}\,\,2\,\,2\,\,2\,\,1\,\,1\,\,\,\,\,\,\,\,\,\,\,$   ${\bf
0}\,\,1\,\,1\,\,1\,\,1\,\,1\,\,\,\,\,\,\,\,\,\,\,$   ${\bf
0}\,\,1\,\,1\,\,1\,\,1\,\,1$

$2\,\,{\bf 0}\,\,0\,\,2\,\,1\,\,1\,\,\,\,\,\,\,\,\,\,\,$   $1\,\,{\bf
0}\,\,0\,\,1\,\,1\,\,1\,\,\,\,\,\,\,\,\,\,\,$   $1\,\,{\bf
0}\,\,0\,\,1\,\,1\,\,1$

$2\,\,0\,\,{\bf 0}\,\,2\,\,1\,\,1\,\,\,\,\,\,\,\,\,\,\,$   $1\,\,
0\,\,{\bf 0}\,\,1\,\,1\,\,1\,\,\,\,\,\,\,\,\,\,\,$
$1\,\,0\,\,{\bf 0}\,\,1\,\,1\,\,1$

$2\,\,2\,\,2\,\,{\bf 0}\,\,1\,\,1\,\,\,\,\,\,\,\,\,\,\,$   
$1\,\,1\,\,1\,\,{\bf
0}\,\,1\,\,1\,\,\,\,\,\,\,\,\,\,\,$
$1\,\,1\,\,1\,\,{\bf 0}\,\,1\,\,1$

$\,1\,\,1\,\,1\,\,1\,\,{\bf 0}\,\,0\,\,\,\,\,\,\,\,\,\,\,$
$1\,\,1\,\,1\,\,1\,\,{\bf 1}\,\,1\,\,\,\,\,\,\,\,\,\,\,\,$
$0\,\,0\,\,0\,\,0\,\,{\bf 0}\,\,0$

$\,1\,\,1\,\,1\,\,1\,\,0\,\,{\bf 0}\,\,\,\,\,\,\,\,\,\,\,$
$1\,\,1\,\,1\,\,1\,\,1\,\,{\bf 1}\,\,\,\,\,\,\,\,\,\,\,\,$
$0\,\,0\,\,0\,\,0\,\,0\,\,{\bf 0}$
\end{center}

Clearly, $d_{ij}+d_{ik}-d_{jk}=p_{ij}+p_{ik}-p_{jk}-p_{ii}=q_{ji}+q_{ik}-q_{jk}$, 
i.e., the triangle inequalities are equivalent on all three  levels:
$d$ - of semimetrics, $p$ - of would-be partial semimetrics and $q$ - of 
would-be quasi-semimetrics.

Now, 
$p_{ij}\ge p_{ii}$ iff $d_{ij}\ge w_i-w_j$ iff 
$q_{ij}\ge 0$; so, (4) is equivalent to (4'),

$p_{ij}\le p_{ii}+p_{jj}$  iff $d_{ij}\le w_i+w_j$ iff $q_{ij} \le w_j$;
so, (5) is equivalent to (5'),

and $2p_{ij}\ge p_{ii}+p_{jj}$ iff $d_{ij}\ge 0$ 
iff $q_{ij}+q_{ji}\ge 0$.
This implies

\begin{lemma}
The following statements hold.

(i) $\wPMET_n=P(\WMET_n)$, $\PMET_n=P(\dWMET_n)$ and 

$\sPMET_n=P(\sWMET_{n})$,

(ii) $\wWQMET_n=Q(\WMET_n)$, $\WQMET_n=Q(\dWMET_n)$ and
 
$\sWQMET_n=Q(\sWMET_{n})$.

\end{lemma}
\vspace{1.3mm}

The metric cone  $\MET_n\in \mathbb{R}^{n \choose
2}$ has a unique orbit of  $3{n
\choose 
3}$ facets $Tr_{ij,k}$. Its symmetry group $\Aut(\MET_n)$
 is $\Sym(n)$   $n\neq4$. 
 The number of 
extreme rays (orbits)
of $\MET_n$ is $3$ ($1$), $7$ ($2$), $25$ ($3$), $296$ ($7$), $55226$
($46$)
 for $3\le n \le 7$. 
 $D(\Ri(\MET_n))=2$ for $n>3$, while $\Ri(\MET_3)=\Sk(\MET_3)=K_3$.  
 $D(\Sk(\MET_n))$ is $1$ for $n=4$,  
$2$ for $5\le
n\le 6$ and  $3$ for $n=7$.

For a partition $\mathcal{S}=\{S_{1}, \dots , S_{t}\}$  of $V_n$,
the  {\em multicut} $\delta(\mathcal{S})\in \MET_n$ has
$\delta_{ij}(\mathcal{S})=1$ if $|\{i,j\}|=2>|\{i,j\}\cap S_h|, 1\le h\le t$ and 
 $\delta_{ij}(\mathcal{S})=0$, otherwise.
Call 
 $\delta(\mathcal{S})$ a {\em $t$-cut} if $S_h\neq\emptyset$ for $1\le h\le t$.
Clearly,  $$\delta(S_1, \dots,
S_t)=\frac{1}{2}\sum_{h=1}^t\delta(S_h,\overline{S_h}).$$

Denote by $\CUT_n$ the cone   generated by all $2^{n-1}-1$ $2$-cuts $\delta(S, 
\overline{S})=\delta(S)$. $\CUT_n=\MET_n$ holds
 for $n\le 4$ and $\Aut(\CUT_n)=\Aut(\MET_n)$.   
The number of facets (orbits) of
$\CUT_n$ is $3$ ($1$), $12$ ($1$), $40$ ($2$), $210$ ($4$), $38780$
($36$)
 for $3\le n\le 7$. $D(\Sk(\CUT_n))=1$ and  $D(\Ri(\CUT_n))=2,3,3$ for
$n=5,6,7$.
See, for example,~\cite{DL,DDF,Du08} for details on  $\MET_n$
 and $\CUT_n$.

The number of $t$-cuts  of $V_n$
 is the number of ways  to partition a set of $n$ objects into
$t$ groups,  i.e., the {\em Stirling number of the
second kind} $S(n,t)=\frac{1}{t!}\sum_{j=0}^t(-1)^j{t \choose j}(t-j)^n$.
So, $S(n,2)=2^{n-1}-1$ and $S(n,n-1)={n\choose 2}$. The number of multicuts  of $V_n$
is the {\em Bell number} $B(n)=\sum_{t=0}^nS(n,t)=\sum_{t=0}^{n-1}(t+1)S(n,t)=\sum_{t=0}^{n-1}
{n-1 \choose t}B(t)$ (the sequence A000110 $=1, 1, 2, 5$, $15, 52, 203, 877, \dots$
in~\cite{Sl}). The number of ways to write $i$ as
a sum of positive integers is  $i$-th {\em partition
number} $Q_i$ (the sequence $A000041$  in~\cite{Sl}).
\vspace{0.7mm}

The $0,1$-valued elements $d\in \MET_n$ are all $B(n)$  multicuts $\delta(\{S_{1}, \dots , 
S_{t}\})$  of $V_n$. It follows by induction  
using
that $d_{1i}=d_{1j}=0$ implies $d_{ij}=0$   and $d_{1i}\neq d_{1j}$ implies $d_{ij}=1$.
In fact, $S_1, \dots, S_t$ are the equivalence classes of the
equivalence $\sim$ on $V_{n}$, defined by $i\sim j$ if $d_{ij}=0$.    

$R(0,1$-$\MET_n)$ consists of all  $S(n,2)$ $2$-cuts; so, $0,1$-$\MET_n=\CUT_n$.

\section{Description of $\wPMET_n$ and $\sPMET_n$}

Denote $\MET_{n;0}=\{(d;(0)): d\in \MET_n\}$ and $\CUT_{n;0}=\{(d;(0)): d\in 
\CUT_n\}$.
So, $\MET_n\simeq \MET_{n;0}\simeq$ $P(\MET_{n;0})\simeq Q(\MET_{n;0})$ and 
$\CUT_n\simeq \CUT_{n;0}\simeq$ $P(\CUT_{n;0})\simeq Q(\CUT_{n;0})$. 
Denote by $\WCUT_n$ the cone $\{d;w)\in \WMET_n: d\in \CUT_{n}\}$ of 
{\em weighted $\ell_1$-semimetrics} on $V_n$.

Denote  $e_j=(((0));w=(w_i=1_{i=j}))\in \WMET_n$.
So, $2P(e_j)=2$ on the position $(jj)$,  $1$ on $(ij),(ji)$ with 
$i\neq j$  and $0$, else; $2Q(e_j)=-1$   on the positions $(ji)$,  $1$ on $(ij)$ (with
$i\neq j$ again) and $0$, else.

\begin{theor}
The following statements hold.

(i) $R(\WMET_n)=\{e_j: j\in V_n\}\cup R(\MET_{n;0})$,

$R(\wPMET_n)=\{2P(e_j): j\in V_n \}\cup P(R(\MET_{n;0}))$,

$R(\wWQMET_n)=\{2Q(e_j): j\in V_n\}\cup Q(R(\MET_{n;0}))$.

(ii) $F(\WMET_n)=\{w_j\ge 0: j\in V_n\}\cup F(\MET_{n;0})$,

$F(\wPMET_n)=\{L_{jj}=p_{jj}\ge 0: j\in V_n\}\cup F(P(\MET_{n;0}))$,

$F(\wWQMET_n)=\{w_j\ge 0: j\in V_n\}\cup F(Q(\MET_{n;0}))$.

(iii) $\Inc(2P(e_j))=|F(\wPMET_n)|-1$ and $\Inc(L_{jj})=|R(\wPMET_n)|-1$.

(iv) $\Ri(\WMET_n)=\Ri(\wPMET_n)=\Ri(\wWQMET_n)=K_n\times \Ri(\MET_n)$, 

$\Sk(\WMET_n)=\Sk(wPMET_n)=\Sk(\wWQMET_n)=K_n\times \Sk(\MET_n)$.

(v) $\wPMET_n$  has  $\Aut$, $D(\Sk)$,  $D(\Ri)$
and edge-connectivity of  $\MET_n$.

(vi)
The  $0,1$-valued elements of $\wPMET_n$ are   
the $B(n+1)$
$0,1$-valued elements of $\PMET_n$ and
 $0,1$-$\wPMET_n=\CUT_{n;0}$.

(vii)
The  $0,1$-valued elements of $\WMET_n$ are
$2^nB(n)$ $0,1$-weighted
multicuts of $V_n$ and $0,1$-$\WMET_n=\WCUT_n$.

$R(\WCUT_n)=\{e_j: j\in V_n\}\cup R(\CUT_{n;0})$ and $\Sk(\WCUT_n)=K_{n+S(n,2)}$.

$F(\WCUT_n)=\{w_j\ge 0: j\in V_n\}\cup F(\CUT_{n;0})$ 
and $\Ri(\WCUT_n)=K_n\times \Ri(\CUT_n)$ has diameter $2$.

\end{theor}

{\em Proof.} 

{\em (i)}.  Let $p\in \wPMET_n$.  We will show that 
$p'=p-\frac{1}{2}\sum_{t=1}^np_{tt}2P(e_t)\in \MET_{n;0}$.  
 (For example, a well-known 
weak
partial
semimetric $i+j$ is the sum of $\sum_tt2P(e_t)$ and the all-zero 
semimetric $((0))$.)
 
In fact, $p'_{ii}=p_{ii}-\frac{1}{2}p_{ii}2P(e_i)_{ii}=0$. Also, 
$p'$ satisfies to 
all triangle inequalities (3), since for different $i,j,k\in V_n$, we have 
$Tr_{ij,k}=p'_{ij}+p'_{ik}-p'_{jk}=$
$$=\left(p_{ij}-\frac{p_{ii}+p_{jj}}{2}\right) 
+\left(p_{ik}-\frac{p_{ii}+p_{kk}}{2}\right)-\left(p_{jk}-\frac{p_{jj}+p_{kk}}{2}\right) = $$
$$ = p_{ij}+p_{ik}-p_{jk}-p_{ii}\geq 0.$$

So,  
$2P(e_i)$, $1\le i\le 
n$,  and the generators of $P(\MET_{n;0})\simeq \MET_{n;0}$
(i.e.,  the $0$-extensions of the generators  of 
$\MET_n$) generate $\wPMET_n$.
They are, moreover, the  generators of $\wPMET_n$ since
 they belongs to 
all  $n$ 
(linearly independent) facets 
$L_{ii}$;
 so, their  rank in $\mathbb{R}^{n+1 \choose
2}$is $({n 
\choose
2}-1)+n={n+1 \choose
2}-1$.

 Clearly, any   $2P(e_i)$ belongs to all facets of 
$\wPMET_n$ except 
 $L_{ii}$, i.e., its incidence is $(n-1)+3{n \choose
3}$. So, its rank  in 
$\mathbb{R}^{n+1 \choose 2}$ is ${n+1 \choose 2}-1$.  
For $\WMET_n$ and $\wWQMET_n$, (i) follows similarly, as well as (ii).

{\em (iii), (iv)}. The ray of $2P(e_i)$ is adjacent to any other extreme
ray $r$, as the set of facets that contain $r$ (with rank
${n+1 \choose 2} - 1$) only loses one element if we intersect
it with the set of facets that contain $2P(e_i)$.

{\em (v)}. The diameters of $\Ri(\wPMET_n)$ and $\Sk(\wPMET_n)$ being $2$,
 their edge-connectivity 
 is equal
to their  minimal degrees~\cite{Ple}. But this degree is the same as of  $\Ri(\MET_n)$ (which 
is 
regular of degree
$\frac{(n-3)(n^2-6)}{2}$
 if $n>3$) and of $\Sk(\MET_n)$, respectively.
$\Aut(\wPMET_n)$ for $n\ge 5$ is $\Sym(n)$, 
because it contains $\Sym(n)$ but cannot be larger than  
$\Aut(\MET_n)=\Sym(n)$. 

{\em (vi)}. If $p\in \wPMET_n$ is $0,1$-valued,
then $p_{ij}=0<p_{ii}=1$ is impossible because   
$2p_{ij}\ge p_{ii}+p_{jj}$; so, $p\in \PMET_n$.
(vii) is implied by (i), (ii).
\qed
\vspace{1.3mm}

Any partial semimetric $p\in \PMET_n$ induces the partial 
order on $V_n$ by defining
 $i\preceq j$ if $p_{ii}=p_{ij}$. This {\em specialization order} is 
important in Computer Science applications, where the  partial metrics 
 act on certain posets called {\em Scott domains}.  
In particular, $i_0\in V_n$ is a {\em $p$-maximal} element in $V_n$ if
$p_{ii}=p_{ii_0}$ for all
$i\neq i_0$. It is a {\em $p$-minimal} element in $V_n$ if
$p_{i_0i_0}=p_{ii_0}$ for all
$i\neq i_0$. 
The  {\em lifting} of $p\in \PMET_n$ is the
function $p^{+}=((p^{+}_{ij}))$, $i,j\in \{0\} \cup V_n$, with
 $p^{+}_{00}=0$, $p^{+}_{0i}=p^{+}_{i0}=p^{+}_{ii}$ for $i\in V_n$ and
$p^{+}_{ij}=p_{ij}$ for
$i,j\in V_n$.  
Clearly,   $0$
is a $p^{+}$-maximal element in the specialization order, induced on 
$\{0\}\cup V_n=\{0,1,\dots, n\}$ by $p^{+}$, since 
$p^{+}_{ii}=p_{ii}$ as well as  $p^{+}_{i0}=p_{ii}$ for all
$i\in V_n$.

\begin{theor}
The following statements hold.

(i) $\sPMET_n=\{p\in \PMET_n: p^{+}\in  \PMET_{n+1}\}$.

(ii)  $\sWMET_n=\MET_{n+1}\simeq P(\MET_{n+1})=\sPMET_n$. 

 (iii)
The $0,1$-valued elements of $\sPMET_n$ are $((0))$ and
$2^{n}-1$  partial $2$-cuts
$\gamma(S\neq \emptyset; \overline{S})$ generating $0,1$-$\sPMET_{n}=CUT_{n+1}$,

$\CUT_{n+1}\simeq P(\CUT_{n+1})=0,1$-$\sPMET_n$ and 
$Q(\CUT_{n+1})=\OCUT_n$.

\end{theor}
{\em Proof.}
We should check for $p^{+}$ only inequalities
(2), (3), (4) involving the new point $0$. $2n+1$ of the required 
inequalities hold as 
equalities: $p^{+}_{00}=0$ and 
$p^{+}_{0i}=p^{+}_{i0}=p^{+}_{ii}=p_{ii}$ for $i\in V_n$.
All $Tr_{0j,i}=p_{ij}-p_{ij}\ge 0$ hold since (4) is satisfied.
All $Tr_{ij,0}=p_{ii}+p_{jj}-p_{ij}\ge 0$ hold 
 whenever  $p$
satisfies (5), i.e, $p\in \PMET_n$.
\qed
\vspace{0.7mm}

Given $p\in \sPMET_n$, the semimetric $P^{-1}(p^{+})\in \MET_{n+2}$ is 
$P^{-1}(p)\in 
\MET_{n+1}$ with the first point split in two coinciding points.
The cone $\sPMET_n$ is nothing but the linear image 
$P(\sWPMET_n=\MET_{n+1})$. So, for $n\ge 4$, $\Aut(\sWPMET_n)=\Sym(n+1)$ on 
$\{0,1,\dots, n\}$ acting as $p'=P(\tau(P^{-1}(p)))$ on $\sPMET_n$ for any 
 $\tau\in \Sym(n+1)$. If $\tau$ fixes  $0$, then
 $p'=\tau(p)$.

\section{$0,1$-valued elements of  $\PMET_n$ and  $\dWMET_n$}

For a partition $\mathcal{S}=\{S_{1}, \dots , S_{t}\}$  of $V_n$ and $A\subseteq \{1, \dots, t\}$,
 let us denote $\hat{A}=\bigcup_{h\in A} S_h$  and $w(\hat{A})=(w_i=1_{i\in \hat{A}})$.
So, weight is constant on each $S_{h}$.

For any $S \subset V_n$, denote  by $J(S)$
the $0,1$-valued function  with
$J(S)_{ij}=1$ exactly
when $i,j \in S_0$. So,  $J(V_n)$ and $J(\emptyset)$ are all-ones and
all-zeros
partial semimetrics, respectively.

For any $S_0 \subset V_n$ and partition $\mathcal{S}=\{S_1,
\ldots, S_{t}\}$ of
$\overline{S_{0}}$,
denote
 $J(S_0)+\delta(S_{0}, S_1, \dots, S_{t})$ by
 $\gamma(S_{0}; S_1, \dots, S_{t})$ and call it
a {\em partial multicut}
or, specifically,
a {\em  partial $t$-cut}.
Clearly, $\gamma(S_{0}; S_1, \dots, S_{t})\in \PMET_n$ and it is
 $P(d;w)$, where $d=2\delta(S_{0}, S_1, \dots, 
S_{t})-\delta(S_{0})=\sum_{i=1}\delta(S_i)$
 and $w=(w_i=1_{i\in S_0})$.

\begin{theor}
The following statements hold.

(i) The $0,1$-valued elements of $\dWMET_n$  are
$\sum_{t=1}^n2^tS(n,t)$ $(\delta(\mathcal{S});w(\hat{A}))$.

$R(0,1$-$\dWMET_n$) consists of all
 such  elements with $(|A|,t-|A|)=(1,0),$ $(0,2)$ or $(1,1)$, i.e., 
$(((0));(1))$ and $2$-cuts
$(\delta(S);w)$ with weight $(0), w'=(1_{i\in S})$ or $w''=(1_{i\notin 
S})$.
There are  $1+3(2^{n-1}-1)$ of them,  in $\lfloor\frac{3n}{2}\rfloor$ 
orbits.

(ii) The $0,1$-valued elements of $\WQMET_n$  are
$Q(2\delta(\mathcal{S})-\delta(\hat{A};w(A))=\delta(\mathcal{S})-\delta'(\hat{A})$.

$R(0,1$-$\dWMET_n$) consists of all such
elements with either $|A|=t-|A|=1$ (o-$2$-cuts), or $2\le
|A|,t-|A|\le n-2$.

(iii) The $0,1$-valued elements of $\PMET_n$  are the  partial multicuts

$P((2\delta(\mathcal{S})-\delta(\hat{A});w(A))=
\delta(\mathcal{S})+J(\hat{A})$ with $|A|\le 1$.

There are $B(n+1)=$ $\sum_{i=0}^n
{n \choose i}B(i)$ of them, in
  $\sum_{i=0}^n Q(i)$ orbits.

$R(0,1$-$\PMET_n$) consists of all such  elements except
$B(n)-(2^{n-1}-1)$, in $Q(n)-\lfloor\frac{n}{2}\rfloor$ orbits,
those   ({\em partial $t$-cuts}) with $|A|=0,t\neq 2$.

\end{theor}

{\em Proof.}

{\em (i)}.
The $0,1$-valued elements of $\WMET_n$ are $0,1$-weighted
multicuts.
Now, the inequality (4') $d_{ij}\ge w_i-w_j$, valid on $\dWMET_n$, implies 
that $w$ is constant on each $S_h$.
$(((0));(1))$ belongs to $R(0,1$-$\dWMET_n)$
since its rank is
${n\choose 2}$ plus $n-1$, the rank of the set of
equalities $d_{ij}= w_i-w_j$.
The same holds for
$(\delta(S);w)$ with weight $(0), w'=(1_{i\in S})$ or $w''=(1_{i\notin 
S})$,
  since their rank is
${n\choose 2}-1$ plus $k\ge 1$
 equalities $w_i=0$ plus, if $k<n$, $n-k$ equalities 
$d_{ij}=w_j-w_{i'}=1$,
 where $w_{i'}=0$ and  $w_{j}=1$.
But the all-ones-weighted $2$-cut $\delta$ is equal to
$\frac{1}{2}((\delta;w')+(\delta;w'')+(((0)));(1))$.
No other
 $(\delta(S_1, \dots, S_t);w)$ belongs to $R(0,1$-$\dWMET_n)$
  since $t$ should be $2$ (otherwise, the rank
will be $<{n\choose 2}-1 +n$) and the weight should be constant on each
$S_h$, $1\le h\le t$.

{\em (ii)}. 
Let $q\in \WQMET$ be  $0,1$-valued.  Without loss of generality, let
$\min_{i=1}^n (w_i)=w_1=0$. But $q_{1i}+w_1= q_{i1}+w_i$ for any $i>1$.
So, $w_i=1$ if and only if $q_{1i}\neq
q_{i1}$. The quasi-semimetrics $q$, restricted on the sets $\{i: w_i=0\}$ 
and
$\{i: w_i=1\}$, should be $0,1$-valued semimetrics, i.e.,  multicuts.

{\em  (iii)} is proven in~\cite{DeDe}.  For example, there are
$52=1\times 1 + 4\times 1 + 6\times 2 + 4\times 5 + 1\times 15$ (1+1+2+3+5 orbits)
$0,1$-valued elements of $\PMET_4$. Among them,
only
$\delta(S_1, \dots, S_{t})$ with $t=1,3,4$, i.e.,
$((0))$, $\delta(\{1\},\{2\},\{3\},\{4\})$ and $6$ elements of the
orbit with $t=3$ are not representatives of extreme rays.   \qed

\section{Two generalized hypermetric  cones}

For a sequence $b=(b_1,\dots,b_n)$ of integers, where $\Sigma_b$ denotes 
$\sum_{i=1}^nb_i$,
and a symmetric $n\times n$ matrix $((a_{ij}))$,
 denote by $H_b(a)$ the sum $-\sum_{1\le i,j\le n}b_ib_ja_{ij}$
 of the entries of the matrix $-b^Tab$.
 
The cone $\HYP_n$ of all {\em hypermetrics}, i.e., semimetrics $d\in \MET_n$
 with  $H_b(d)\ge 0$, whenever $\Sigma_b=1$, was introduced in~\cite{De60}.
 
This cone  is polyhedral~\cite{DGL93}; $\HYP_n\subseteq \MET_n$ with 
equality for $n\le 
4$ and $\CUT_n\subseteq \HYP_n$  with equality
 for $n\le 6$.
$\HYP_7$  was described in~\cite{DeDu04}.
 
The hypermetrics have deep connections with Geometry of
Numbers and Analysis; see, for example,~\cite{DeTe87,DeGr93,DGL95} and
Chapters 13-17, 28 in~\cite{DL}. So,  generalizations of $\HYP_n$ can put  
 those connections in a more general setting.

For a weighted semimetric $(d;w)\in \WMET_n$, we will use the notation:

$\Hyp_b(d;w)=\frac{1}{2}H_b(d)+(1-\Sigma_b)\langle b,w\rangle \ge 0$ and

$\Hyp'_b(d;w)=\frac{1}{2}H_b(d)+(1+\Sigma_b)\langle b,w\rangle \ge 0$.
\vspace{0.7mm}

Denote by $\WHYP_n$ the cone of all {\em weighted hypermetrics}, i.e., $(d;w)\in \WMET_n$
 with  $\Hyp_b(d;w)\ge 0$ and $\Hyp'_b(d;w)\ge 0$ for all $b$ with $\Sigma_b = 1$ or $0$.
Denote by $\PHYP_n$ the cone of all {\em partial hypermetrics}, i.e., $p\in \wPMET_n$ with 
 $\Hyp_b(P^{-1}(p))\ge 0$ for all $b$ with $\Sigma_b = 1$ or $0$. For $p=P(d;w)$ and
$(q,w)=Q(d;w)$, we have
$$\Hyp_b(d;w)=H_b(p)+\sum_{i=1}^nb_ip_{ii}=
H_b(q)+(1-\Sigma_b)\langle b,w\rangle .$$
$\WHYP_n\subset \dWMET_n$ and $\PMET_n\supset \PHYP_n$ hold  since the 
needed
inequalities
$w_{i}\ge 0$, (4') (and  (4)) are provided by permutations of $\Hyp'_{(1, 0,\dots, 0)}(d;w)\ge
0$ 
 and $\Hyp_{(1, -1, 0,\dots, 0)}(d;w)\ge 0$.

\begin{lemma} 

Besides 
 the cases 
 $\PMET_3=\PHYP_3=0,1$-$\PMET_3$ and  $0,1$-$\dWMET_n= 
\WHYP_n$ for $n=3,4$, $0,1$-$\dWMET_n\subset 
\WHYP_n\subset \dWMET_n$
$\simeq$ 
 $\PMET_n\supset \PHYP_n\supset 0,1$-$\PMET_n$ holds.

\end{lemma}

{\em Proof.}

Denoting $\langle b, (1_{i\in S_h})\rangle$ by $r_h$, we have
$r_0=\Sigma_b-\sum_{h=1}^tr_h$ and
 $$H_b(\delta(S_0,S_1,\dots,
S_t))=\frac{1}{2}\sum_{h=0}^tH_b(\delta(S_h,\overline{S_h}))=
\sum_{h=0}^tr_h(r_h-\Sigma_b).$$

Let $(d=2\delta(S_0,S_1,\dots, S_t)-\delta(S_0); w=(1_{i\in S_0}))$ be a
generic $P^{-1}(p)$, where $p$ is $0,1$-valued element of $\PMET$
belonging to its extreme ray. Then
$\frac{1}{2}H_b(d)=\frac{1}{2}(2\sum_{h=0}^tr_h(r_h-\Sigma_b)-
2r_0(r_0-\Sigma_b)=\sum_{h=1}^tr_h(r_h-\Sigma_b)$ implies

$\Hyp_b(d;w)=\sum_{h=1}^tr_h(r_h-1) -\Sigma_b(\Sigma_b-1)\ge 0$ for our $\Sigma=0,1$.

All $0,1$-valued elements $(d;w)$ of $\dWMET_n$ belonging to its extreme
rays are
$(((0)); (1))$, $(\delta(S); (0))$,  $(\delta(S,); w'=(1_{i\in S}))$ and
$(\delta(S); w''=(1_{i\notin S}))$. For them,
$\Hyp'_b(d;w)=(\Sigma_b+1)\Sigma_b$,
$r_S(r_S-\Sigma)$, $r_S(r_S+1)$ and $(\Sigma_b-r_S)(\Sigma_b-r_S+1)$ hold, respectively, 
and so, $\Hyp'_b(d;w)\ge 0$ for our $\Sigma=0,1$. 

Assuming polyhedrality of $\WHYP_n$, the cases $n=3,4$ were checked by
computation; see Lemma below. \qed

\begin{lemma} The following statements hold.

(i) All  facets of $\WHYP_n$, $n\le 4$,
up to $\Sym(n)$ and $0$-extensions, are
$\Hyp_b$ with $b=(1,-1), (1,1,-1),(1,1,-1,-1)$ and $\Hyp'_b$ with 
$b=(1),$ $(1,1,-1)$, $(1,1,1,-2)$, $(2,1,-1,-1)$.

(ii) Besides $w_i\ge 0$, among the facets of $P^{-1}(\PHYP_n)$, $n\le 5$,
up to $\Sym(n)$ and $0$-extensions, are:
$\Hyp_b$ with $b=(1,-1)$, $(1,1,-1)$, $(1,1,-1,-1)$, $(1,1,1,-1,-1)$, $(1,1,1,-1,-2)$,
$(2,1,1,-1,-1)$.
\end{lemma}

{\em Proof.}

It was obtained by direct computation.
The equality $\WHYP_n =0,1$-$\WMET_n$ for $n=3,4$ holds, because only inequalities which are
 requested in  $\WHYP_n$ appeared among those of $0,1$-$\WMET_n$.

The facets of $\PHYP_4$ were deduced by computation using the tightness of the 
inclusions $0,1$-$\PMET_4\subset
\PHYP_4\subset \PMET_4$ (see  
Table~\ref{tab:tablL1-4facets}):
$0,1$-$\PMET_4$
contained exactly one facet (orbit $F_5$) different from $\Hyp_b$ and $p_{ii}\ge0$, and
$\PMET_4$
contained exactly two (orbits $R_{10}$ and $R_{11}$) non-$0,1$-valued extreme ray 
representatives.
The $6$ rays from $R_{10}$ are removed by $6$ respective $\Hyp_b$ with 
$b=(1,1,-1,-1)$, while the $12$ rays   from
$R_{11}$ are removed by $12$  $F_5$.
\qed
\vspace{0.7mm}

\begin{table}
\begin{center}
\scriptsize
\makebox[\textwidth]{
\begin{tabular}{|c|c|cccccccccc|c|c|c|}
 \hline
$R_i$ & Representative $p$
&11&21&22&31&32&33&41&42&43&44&Inc.&Adj.&$|O_i|$\\
 \hline
$R_1$&$\gamma(\{1,2,3,4\};)$&1&1&1&1&1&1&1&1&1&1&24&20&1\\
$R_2$&$\gamma(\{2\};\overline{\{2\}})$  &0&1&1&0&1&0&0&1&0&0&21&38&4\\
$R_3$&$\gamma(\overline{\{3\}}; \{3\})$ &1&1&1&1&1&0&1&1&1&1&19&17&4\\
$R_4$&$\gamma(\emptyset;\overline{\{3\}}, \{3\})$
&0&0&0&1&1&0&0&0&1&0&19&32&4\\
$R_5$&$\gamma(\{1,2\};\overline{\{1,2\}}) $ &1&1&1&1&1&0&1&1&0&0&18&31&6\\
$R_6$&$\gamma(\emptyset; \{1,2\},\overline{\{1,2\}}) $ &0&0&0&1&1&0&1&1&0&0&16&32&3\\
$R_7$&$\gamma(\{1,4\};\{2\}, \{3\}) $ &1&1&0&1&1&0&1&1&1&1&14&14&6\\
$R_8$&$\gamma(\{1\};\{2\}, \{3,4\}) $&1&1&0&1&1&0&1&1&0&0&14&20&12\\
$R_9$&$\gamma(\{4\};\{1\},\{2\},\{3\}) $ &0&1&0&1&1&0&1&1&1&1&9&9&4\\
 \hline
$R_{10}$& &1&1&0&1&1&0&2&1&1&1&10&18&6\\
 \hline
$R_{11}$& &0&2&0&1&1&0&2&2&3&2&9&9&12\\
 \hline 
 \hline

$F_1$&$L_{11}: p_{11}\ge 0$& 1&0&0&0&0&0&0&0&0&0&29&36&4\\
$F_2$&$\Hyp_{(-1,1,1,0)}
\ge 0$&
 -1&1&0&1&-1&0&0&0&0&0&26&24&12\\
$F_3$&$M_{12}=\Hyp_{(-1,1,0,0)}  
\ge 0$& -1&1&0&0&0&0&0&0&0&0&23&23&12\\
\hline
$F_4$&$\Hyp_{(1,1,-1,-1)}\ge 0$& 0&-1&0&1&1&-1
&1&1&-1&-1&16&12&6\\
\hline
$F_5$&$
H_{(2,1,-1,-1)}-2p_{11}\ge 0$ & -2&-2&0&2&1&0
&2&1&-1&0&9&9&12\\ \hline 
\hline
\end{tabular}
}
\caption{The 
  orbits of extreme rays in $\PMET_{4}$
and facets in $0,1$-$\PMET_{4}$}
\label{tab:tablL1-4facets}
\end{center}
\vspace{-5mm}
\end{table}

\section{Oriented multicuts and quasi-semimetrics}

For an ordered partition  $(S_{1}, \dots , S_{t})$
of $V_n$ into non-empty subsets, the {{\em oriented multicut}}
(or {{\em o-multicut}}, {\em o-$t$-cut})  $\delta^{'}(S_{1},\dots,S_{t})$ on
$V_n$ is defined by:

$\delta^{'}_{ij}(S_{1}, \dots , S_{t})=
\left\{
\begin{array}{ccc} 1, & \mbox{if} &
i \in S_{h}, j \in S_{m}, m>h, \\ 0, && \mbox{otherwise}. \end{array}\right.$

The o-$2$-multicuts
$\delta^{'}(S,\overline{S})$ are  called {\em o-cuts} and
denoted by $\delta^{'}(S)$. Clearly, $$\delta (S_{1}, \dots,
S_{t})=\sum_{i=1}^t \delta^{'}(S_{i})=\sum_{i=1}^t \delta^{'}
(\overline{S_{i}})=\frac{1}{2}\sum_{i=1}^t \delta (S_{i}).$$ 
Denote by $\OCUT_n$ and $\OMCUT_n$
 the cones generated by  $2^{n}-2$  non-zero o-cuts and $Bo(n)-1$
non-zero o-multicuts, respectively. (Here,
 $Bo(n)$ are the {\em ordered Bell numbers} given by the sequence
$A000670$
in~\cite{Sl}.) So, $\CUT_n=\{q+q^T:q\in \OMCUT_n\}$. In general, $Z_2\times
\Sym(n)$ is a symmetry group of $\QMET_n$, $\OMCUT_n$, $\WQMET_n$, $\OCUT_n$;
Dutour, 2002, proved that it is the full  group of those cones.
The cones $\QMET_n$ and  $\OMCUT_n$ were studied in~\cite{DDP}.  
Clearly, $\delta^{'}_{ij}(S_{1}, \dots , S_{t})\in \WQMET_n$    
if and only if
$t=2$ and then $w=(1_{i\notin S_1})$. So, $\OCUT_n=\OMCUT_n\cap \WQMET_n$.

\begin{theor}
The following statement holds.

$\OCUT_n=Q(\CUT_{n+1}=0,1$-$\sWMET_n)=Q(0,1$-$\dWMET_n)=$
$0,1$-$Q(\dWMET_n)$.

\end{theor}

{\em Proof.}
Given a representative  $(d;w)=(\delta(S);w')$ ,
$(\delta(S);w'')$, $(\delta(S);(0))$, $(\delta(\emptyset);(1))$ of an extreme ray of
 $ 0,1$-$\dWMET_n$, we have $Q(q;w)=(\delta'(S),w')$,
$(\delta'(\overline{S}),w'')$,
$(\delta(S),(0))$, $(\delta(\emptyset),(1))$, respectively.
But $\delta(S)=\delta'(S)+\delta'(\overline{S})$  and $(((0)),t(1))$ are
not extreme rays.\qed
\vspace{0.7mm}

The above equality $\OCUT_n=Q(0,1$-$\sWMET_n)$ means that $q\in \OCUT_n$ 
are
$Q(d;w)$, where $(d,w)$ is a semimetric $d'\in \CUT_{n+1}$ on
$V_{n}\cup\{0\}$.
 So, $q_{ij}=\frac{1}{2}(d_{ij}'-d_{0i}'+d_{0j}')$.
 But $\CUT_{n}$ is the set of {\em $\ell_1$-semimetrics} on $V_{n}$, see~%
\cite{DL}.
So, $q\in \OCUT_n$ can be seen as {\em $\ell_1$-quasi-semimetrics}; it was
realized in
\cite{DDP,CMM06}. In fact, $\OCUT_n$ is the set of
quasi-semimetrics $q$ on $V_n$, for which there is  some $x_1, \dots,
x_m\in \mathbb{R}^m$  with all $q_{ij}=||x_i-x_j||_{1.or}$, where the
{\em oriented
$\ell_1$-norm} is defined as
$||x-y||_{1.or}=\sum_{k=1}^m\max (x_k-y_k,0)$; the proof is the same as  in
Proposition 4.2.2 of~\cite{DL}.
\vspace{1.3mm}

Let $C$ be any cone closed under {\em reversal}, i.e., $q\in C$ implies
 $q^T\in C$.
If the linear inequality $\sum_{1\le i,j\le n}f_{ij}q_{ij}=\langle F,q \rangle 
\ge 0$ 
is valid on 
$C$, then 
$F$ also defines a face of $\{q+q^T: q\in C\}$.
Given a valid inequality  $G: \sum_{1\le i<j\le n}g_{ij}d_{ij}$ of 
$\{q+q^T: q\in C\}$ and  an {\em oriented $K_n$} (i.e.,  exactly  one arc connects any 
$i$ and $j$) 
$O$, let $G^O=((g_{ij}^O))$ where  
 $g_{ij}^O=g_{ij}$ if 
the arc $(ij)$ belongs to $O$ and $=0$, otherwise.
Call $G^O$ {\em standard} if there exists  $\tau \in \Sym(n)$ 
with $(ij) \in O$  if and only if $\tau(i) < \tau(j)$, and  
  {\em reversal-stable} ({\em rs} for short) if   $\langle G^O,q \rangle =\langle 
G^O,q^{T} \rangle$.
In general, $G^O$ is not valid on 
$C$ and does not preserve the rank of $G$.

For example, the standard  
$Tr_{12,3}: q_{13}+q_{23}-q_{12}\ge 
0$ is 
not valid on $\OCUT_n$, and the standard  $L_{ij}: q_{ij}\ge 0$ defines a facet in  
$\OCUT_n$, while $G: d_{ij}\ge 0$ only defines a face in $MET_n$.
If $F=G^O$ is rs,
 then $\langle F,q \rangle =\frac{1}{2}\langle 
G,q+q^{T} 
\rangle$, i.e., $F$ is valid 
on $C$ if $G$ is valid on $\{q+q^T: q\in C\}$.

Let $E$ be an equality that holds on $C$, i.e.,
$\sum_{1 \le i,j \le n} e_{ij} q_{ij} = \langle E, q \rangle = 0$
holds for any $q \in C$. If the dimension of the subspace $\mathcal{E}$,
spanned by all such equalities, is greater than zero, and $F \ge 0$ is a
facet-defining inequality, then for any $E \in \mathcal{E}$, $F+E \ge 0$
defines the same facet.
We call a facet standard or rs if one of its defining inequalities is
standard or rs. Of all the defining inequalities we can choose one of them
(up to a positive factor) to be the {\em canonical representative} -- let it
be such a $G = F+E$, $E \in \mathcal{E}$ that $\langle G, E \rangle = 0$
holds for all $E \in \mathcal{E}$, i.e., $G$ is orthogonal to $\mathcal{E}$.

\begin{lemma}
Let $C$ be a cone closed under reversal, and
$\mathcal{E}$ the subspace of its equalities. Then, the following statements
hold.

(i) If $C \subseteq \WQMET_n$ is of the same dimension as $\WQMET_n$,
then $\mathcal{E}$ is spanned by the equalities
$q_{ij} + q_{jk} + q_{ki} = q_{ji} + q_{kj} + q_{ik}$ for $i,j,k \in V_n$
and its dimension is ${n-1 \choose 2}$.

(ii) For each $E \in \mathcal{E}$, $E$ is rs.

(iii) If a facet of $C$ is rs, then all of its defining inequalities are rs.

(iv) A facet is rs iff its canonical representative $G$ is {\em symmetric},
e.g. $G = G^T$ holds.
\end{lemma}

\clearpage

{\em Proof.}

{\em (i)}. The equalities
$E_{ijk} = q_{ij} + q_{jk} + q_{ki} - q_{ji} - q_{kj} - q_{ik} = 0$ for
$i,j,k \in V_n$ follow directly from the weightability condition
$q_{ij} + w_i = q_{ji} + w_j$. Since $E_{ijk} = E_{jki} = -E_{kji}$ and
$E_{jk\ell} = E_{ijk} - E_{ij\ell} + E_{ik\ell}$ hold, we can choose a basis
of $\mathcal{E}$ such that the indices $(ijk)$ of the basis elements
$E_{ijk}$ are ordered triples that all contain a fixed element of $V_n$
(say, $n$). There are ${n-1 \choose 2}$ such triples, and since all such
$E_{ijk}$ are linearly independent, the subspace $\mathcal{E}$ has dimension
${n-1 \choose 2}$.

{\em (ii)}. As $C$ is closed under reversal, each equality
$E \in \mathcal{E}$ holds for both $q, q^T \in C$. Therefore,
$0 = \langle E, q \rangle = \langle E, q^T \rangle$, so $E$ is rs.

{\em (iii)}. If $F$ is a defining inequality of a facet and is rs, then
for each $E \in \mathcal{E}$ and $q \in C$,
$\langle F+E, q^T \rangle = \langle F+E, q \rangle$,
so $F+E$, and by extension any defining inequality, is also rs.

{\em (iv)}. Clearly, if $G$ is symmetric, it is also rs and so is the facet
it defines. If a facet is rs, then by (iii), so is its canonical
representative $G$, for which
$\langle G, q \rangle = \langle G, q^T \rangle = \langle G^T, q \rangle$
holds for all $q \in C$. Therefore, $G - G^T \in \mathcal{E}$, but as
$G - G^T$ is also orthogonal to $\mathcal{E}$, $G = G^T$ follows.
\qed
\vspace{1mm}

The facets $OTr_{ij,k}$
 and $L_{ij}$
 (only $1$st is rs)
 of $\WQMET_n$ are standard and of the form $\Hyp_b$ where  $b$ is a permutation of
$(1,1,-1,0,\dots, 0)$
or $(1,-1,0,\dots,0)$. 
$\OCUT_4$ has one more orbit:
 six standard, non-rs facets of the form $\Hyp_b$ where  $b$ is a permutation of
$(1,1,-1,-1)$, or $q_{13}+q_{14}+q_{23}+q_{24}-(q_{12}+q_{34})\ge 0$.

$\OCUT_5$ has,  up to $\Sym(n)$, $3$ {\em new} (i.e., in addition to  
$0$-extensions of the 
facets of $\OCUT_4$) such orbits: one standard rs  
 $(1,1,1,-1,-1)$  and two non-standard, non-rs orbits. 
$\OCUT_6$ has, among its $56$ new orbits, 
 two  non-standard rs orbits for 
$b=$ $(2,1,1,-1,-1,-1)$ and  $(1,1,1,1,-1,-2)$.

The adjacencies of cuts in $\CUT_n$ 
 are defined only by the 
facets $Tr_{ij,k}$, and adjacencies of those facets are  defined only by cuts.
It gives at once  ${n\choose 2}-1$ linearly independent facets $OTr_{ij,k}$ containing
any given pair $(\delta^{'}(S_1),\delta^{'}(S_2))$, using that $OTr_{ij,k}$ are rs 
facets. So, only $n$ more facets are needed to get the adjacencies of o-cuts. It is a 
way to prove Conjecture 1 (i) below.
 \vspace{1mm}

Call a {\em tournament} ($K_n$ with unique arc between any $i,j$)  {\em 
admissible} if its 
arcs can be partitioned  into arc-disjoint directed cycles. 
It does not exists for even $n$, because then the number of arcs 
involving each vertex is 
odd, while each cycle provides $0$ or $2$ such arcs. But for odd $n$, 
there 
are at least $2^{\frac{n-3}{2}}$ admissible tournaments: take the 
decomposition of $K_n$ into $\frac{n-1}{2}$  disjoint  Hamiltonian
cycles and, fixing the order on one them, all possible orders on remaining 
cycles. For odd $n$, denote by $Oc$ the  {\em canonic admissible tournament}  
consiting 
of all $(i,i+k)$ with $1\le i\le n-1, 1\le k\le {\lceil \frac{n}{2}\rceil}+1-i$
and  $(i+k,i)$ with $1\le i\le {\lfloor \frac{n}{2} \rfloor}, {\lceil 
\frac{n}{2}\rceil} \le k\le n-i$, i.e., 
$0=C_{1,2,3,4,5,6,7,\dots}+C_{1,3,5,7,\dots}+C_{1,4,7,\dots}+\dots$.
 The  {\em Kelly conjecture}  
state that the
arcs of every
{\em regular} (i.e., the vertices have the same outdegree)
tournament can be partitioned  into arc-disjoint directed Hamiltonian
cycles.

$0$-extensions of $q_{ij}\ge 0$ and
$q_{13}+q_{14}+q_{23}+q_{24}-(q_{12}+q_{34})\ge 0$
 can be seen, as the first instances (for $b=(1,-1,0,\dots, 0)$, $(1,1,-1,-1,0,\dots, 0)$)
of the {\em
oriented
 negative type inequality}
$\ONeg_{b,O}(q)=-\sum_{1\le i<j\le n}b_ib_jq_{a(ij)}\ge 0$, where 
for a given $b=(b_1,\dots,b_n)\in \mathbb{Z}^n$, $\Sigma_b=0$, and the arcs $a(ij)$ on the 
edges $(ij)$ by some rule defined by a given tournament $O$. 

Denote by $\OWHYP_n$ the cone consisting of all $q\in \WQMET_n$, satisfying the two above
orbits  and all {\em 
oriented
hypermetric inequalities} 
$$\OHyp_{b,O}(q)=-\!\!\!\!\sum_{1\le i<j\le n}b_ib_jq_{a(ij)}\ge 0,$$
where $b=(b_1,\dots,b_n)\in \mathbb{Z}^n$, $\Sigma_b=1$, $O$ is an admissible tournament, and the 
arc $a(ij)$ on the edge $(ij)$
is the same 
as in $O$ if 
$b_ib_j\ge 0$,
or the opposite one otherwise. So, $\OWHYP_n=\OCUT_n$ for $n=3,4$.

\begin{theor}
$\OCUT_n\subset \OWHYP_n\subset \WQMET_n$ holds for $n\ge 5$.

\end{theor}

{\em Proof.}

Without loss of generality, let $b_i=1$ for $1\le i\le \lfloor \frac{n}{2} 
\rfloor$ and $b_i=-1$, otherwise. The general case means only that we have 
sets  of $|b_i|$ coinciding points.
$\OHyp_{b,O}$ is rs,
 because 
Lemma 4 in~\cite{DeDe} implies that any inequality on a $q\in \WQMET_n$
is preserved by the reversal of $q$.
So, $\OHyp_{b,O}(q)=\frac{1}{2}\Hyp_{b}(q+q^{T})$.
On an o-cut  $\delta^{'}(S)$ it gives, putting $r=\langle b,(1_{i\in S}) \rangle$, 
$$\frac{1}{2}\Hyp_{b}(\delta^{'}(S)+\delta^{'}(\overline{S}))=\Hyp_{b}(\delta(S))=
 r(r-\Sigma_b)\ge 0.$$ 
\vspace{0.4mm}

$\OWHYP_5$  has, besides o-cuts,
$40$ extreme rays in two orbits: $F_{ab}, F^{'}_{ab}$, having $2$ on the position $(ab)$, $1$ on
$ba$, $0$ on three other $(ka)$ for $k \ne b$ in $F_{ab}$,
or
on three other $(bk)$ for $k \ne a$ in $F^{'}_{ab}$,
and ones on other non-diagonal positions.
Also, $D(\Sk(\OWHYP_5))=D(\Ri(\OWHYP_5))=2$.
\vspace{0.7mm}

The cone $\QHYP_n=\{q\in \QMET_n: ((q_{ij}+q_{ji}))\in \HYP_n\} $
was considered in~%
\cite{DDP}. Clearly, it is polyhedral and coincides with $\QMET_n$ for 
$n=3,4$;
$\QHYP_5$ has
$90$ facets ($20+60$ from $\QMET_5$ and those with $b=(1,1,1,-1,-1)$) and
$78810$ extreme
rays; $D(\Ri(\QHYP_5))=2$.
\vspace{1.3mm}

Besides $\OCUT_3=0,1$-$\WQMET_3=\WQMET_3$ and $0,1$-$\WQMET_4=\WQMET_4$,
$\OCUT_n\subset 0,1$-$\WQMET_n\subset \WQMET_n$ holds. We conjecture   
 $\Sk(\OCUT_n)\subset \Sk(0,1$-$\WQMET_n)\subset \Sk(\WQMET_n)$ and
$\Ri(0,1$-$\WQMET_n)$ $\supset
\Ri(\WQMET_n)\supset \Ri(\MET_n)$. 
$0,1$-$\WQMET_5$ has
$OTr_{ij,k}$, $L_{ij}$ and $3$  other, all standard, orbits. Those
facets give, for permutations of $b=(1,-1,1,-1,1)$, the
non-negativity of
$-\sum_{1\le i<j\le 5}b_ib_jq_{ij}$ plus $q_{24},q_{23}$ or
$q_{12}+q_{45}$.

The cone $\{q+q^T: q\in 0,1$-$\WQMET_n\}$ coincides with $\MET_n$ for $n\le 5$, 
but for $n=6$ it has $7$ orbits of extreme rays 
 (all those of $\MET_6$ except the one, good 
representatives of which are  not $0,1,2$-valued as required);
 its skeleton, excluding another orbit of $90$ 
rays, is an induced subgraph of   $\Sk(\MET_6)$. 
It has 
$3$ orbits of facets including $Tr_{ij,k}$ (forming $\Ri(\MET_6)$ in its ridge graph) and the orbit of
$\sum_{(ij)\in C_{123456}}d_{ij}+ d_{14}+d_{35}-d_{13}-d_{46}-2d_{25}\ge 0$.
\vspace{0.8mm}

If $q\in \QMET_n$ is  $0,1$-valued
 with  $S=\{i:q_{i1}=1\}$, $S'=\{i:q_{1i}=1\}$, then $q_{ij}=0$ for $i,j\in   
\overline{S}\cap \overline{S'}$ (since $q_{i1}+q_{1j}\ge q_{ij}$)
and $q_{ij}=q_{ji}=1$ for $i\in S, j\in \overline{S'}$ (since
$q_{ij}+q_{j1}\ge
q_{i1}$); so, $|\overline{S}\cap \overline{S'}|(|\overline{S}\cap
\overline{S'}|-1)+
|S|(|\overline{S}|-1)+|S'|(|\overline{S'}|-1)-|S\cap\overline{S'}|
|\overline{S}\cap S'|$
elements $q_{ij}$ with $2\le i \neq j\le n$ are defined.

\section{The cases of $3,4,5,6$ points}

In Table~\ref{tab:MainLovelyTable} we summarize the most important
numeric information on cones under consideration for $n\le 6$. The column $2$
indicates the dimension of the cone, the columns $3$ and $4$ give
the number of extreme rays and facets, respectively; in
parentheses are given the numbers of their orbits. 
The columns $5$ and $6$
give the diameters of the skeleton  and  the ridge graph.
The expanded version of the data can be found on the third author's  
homepage~\cite{Vi}.  

\begin{table}[t]
\begin{center}
\scriptsize
\makebox[\textwidth]{
\begin{tabular}{|c|c|c|c|c|c|}
 \hline \hline
 cone &dim. &Nr.  ext. rays (orbits) &Nr.  facets (orbits)&diam.&diam.  
dual
\\
 \hline
$\wPMET_3$ &6&6 (2)&6 (2) &1&1\\
$\wPMET_4$ &10&11 (3)&16 (2)&1&2 \\
 $\wPMET_5$ &15&30 (4)&35 (2) &2&2\\
$\wPMET_6$ &21&302 (8)&66 (2) &2&2\\
\hline \hline
$\sPMET_3=0,1$-$\sPMET_3$&6&7 (2)&12 (1)&1&2\\
  $\sPMET_4
$&10&25 (3)&30 (1)&2&2 \\
 $\sPMET_5
$ &15&296 (7)&60 (1)&2&2\\
 $\sPMET_6
$ &21&55226 (46)&105 (1)&3&2 \\
\hline
  $0,1$-$\sPMET_4
$&10&15 (2)&40 (2)&1&2 \\
 $0,1$-$\sPMET_5
$ &15&31 (3)&210 (4)&1&3\\
 $0,1$-$\sPMET_6
$ &21&63 (3)&38780 (36) &1&3\\ 
\hline \hline
 $\PMET_3=0,1$-$\PMET_3$ &6&13 (5)&12 (3)&3&2 \\
 $\PMET_4$
&10&62 (11)&28 (3)&3&2 \\
  $\PMET_5
$ &15&1696 (44)&55 (3)&3&2\\
  $\PMET_6$ &21&337092 (734)&96 (3)&3&2\\
\hline
 \hline
$\PHYP_4$&10&56 (10)&34 (4)&3&2\\
\hline
  $0,1$-$\PMET_4$ &10&44 (9)&46 (5) &3&2\\
  $0,1$-$\PMET_5$ &15&166 (14)&585 (15)&3&3\\
  $0,1$-$\PMET_6$ &21&705 (23)&&3&\\
\hline
$0,1$-$\dWMET_3$ &6&10 (4)&15 (4) &2&2\\
$0,1$-$\dWMET_4$ &10&22 (6)&62 (7)&2&3 \\
 $0,1$-$\dWMET_5$ &15&46 (7)& 1165 (27) &2&3\\
$0,1$-$\dWMET_6$ &21& 94 (9)& 369401 (806) &2&\\
\hline \hline
$\WQMET_3=\OCUT_3$ &5&6 (2)&9 (2) &1&2\\
$\WQMET_4=0,1$-$\WQMET_4$ &9&20 (4)& 24 (2)&2&2 \\   
 $\WQMET_5$ &14&190 (11)& 50 (2) &2&2\\
$\WQMET_6$ &20&18502  (77)&90  (2) &&2\\
\hline
 $0,1$-$\WQMET_5$ &14&110  (8)& 250 (5) &2&2\\
$0,1$-$\WQMET_6$ &20&802  (17)& &&\\
\hline
$\{q+q^T: q\in 0,1$-$\WQMET_6\}$ &15&206  (7)& 510 (3) &2&3\\
\hline
$\OWHYP_5$ &14&70 (6)&90 (4)&2&2 \\
\hline
$\OCUT_4$ &9&14 (3)&30 (3)&1&2 \\
 $\OCUT_5$ &14&30 (4)& 130 (6) &1&3\\  
$\OCUT_6$ &20& 62 (5)& 16460 (62) &1&\\
\hline \hline
\end{tabular}
}
\caption{Main parameters of cones with $n\le 6$}
\label{tab:MainLovelyTable}
\end{center}
\vspace{-5mm}
\end{table}

In the simplest case $n=3$
 the numbers of extreme rays 
and facets are:

$0,1$-$\WMET_3=\WHYP_3=\WMET_3\simeq \wPMET_3$: ($6, 6$,  simplicial) and
                   $0,1$-$\wPMET_3$: ($3,3$,  simplicial);

$0,1$-$\sWMET_3=\sWMET_3=\CUT_4=\HYP_4=\MET_4\simeq 0,1$-$\sPMET_3=\sPMET_3$:  
($7, 12$);

$0,1$-$\PMET_3=\PHYP_3=\PMET_3$ ($13, 12$) and $0,1$-$\dWMET_{3}$: ($10, 15$);

$0,1$-$\QMET_{3}=\QHYP_3=\QMET_{3}$: ($12, 12$, simplicial) and 
 $\OCUT_3=0,1$-$\WQMET_{3}=\WQMET_{3}$: ($6, 9$).

$R(\dWMET_3)\setminus R(0,1$-$\dWMET_3)$ and 
$F(0,1$-$\dWMET_3)\setminus 
F(\dWMET_3)$ consist of $3$ simplicial elements forming $\overline{K_3}$ 
in the graph. But only $\Ri(0,1$-$\dWMET_3)$ is an induced subgraph of 
$\Ri(\dWMET_3)$.

Recall that $2^{n-1}-1$ is the Stirling number $S(n,2)$, $\Sk(\CUT_n)=K_{S(n,2)}$, and~\cite{DD}
  $\Ri(\MET_n)$, $n\ge 4$, has diameter $2$
with $Tr_{ij,k}\nsim Tr_{i'j',k'}$ whenever
they are {\em conflicting}, i.e., have values of different sign on a position $(p,q)$,
$p,q\in {\{i,j,k\}\cap \{i',j',k'\}}$. Clearly, $|\{i,j,k\}\cap \{i',j',k'\}|$ should
be $3$ or $2$, and $Tr_{ij,k}$ conflicts with $2$ and $4(n-3)$ $Tr_{i'j',k'}$'s, respectively.
The proofs of the conjectures below should be tedious but easy.

\clearpage

\begin{table}[t]\begin{center}
\scriptsize
\makebox[\textwidth]{
\begin{tabular}{|c|c|cccccc|c|c|c|}
 \hline
$R_i$ & Representative
&11&21&22&31&32&33&Inc.&Adj.&$|R_i|$\\
 \hline
$R_1$ $\blacktriangle$&$\gamma(\{1,2,3\};)$
&1&1&1&1&1&1&9&6&1\\
$R_2$ $\circ$&$\gamma(\{1\};\{2,3\})$ &1&1&0&1&0&0&8&9&3\\
$R_3$ $\bullet$&$\gamma(\{2,3\};\{1\})$ &0&1&1&1&1&1&7&6&3\\
$R_4$ $\square$&$\gamma(\emptyset; \{1\},\{2,3\})$
&0&1&0&1&0&0&7&8&3\\
$R_5$ $\blacksquare$&$\gamma(\{3\};\{1\},\{2\})$
&0&1&0&1&1&1&5&5&3 \\ \hline \hline

$F_1$ $\circ$&$L_{11}: p_{11}\ge 0$& 1&0&0&0&0&0&8&9&3\\
$F_2$ $\blacktriangle$&$Tr_{12,3}:  p_{13}+p_{23}-p_{12}-p_{33}\ge 0$&0&-1&0&1&1&-1&8&7&3\\
$F_3$ $\bullet$&$M_{12}:  p_{12}-p_{11}\ge 0$& -1&1&0&0&0&0&7&6&6\\
 \hline \hline
\end{tabular}
}
\caption{The 
 orbits of extreme rays and facets in
$\PMET_{3}=0,1$-$\PMET_{3}$}
\label{tab:tablL1-3facets}
\end{center}\end{table}
\vspace{0.7mm}

\begin{figure}[t!]
	\centering
    \makebox[\textwidth]{
        \begin{minipage}{0.6\textwidth}
            \includefig{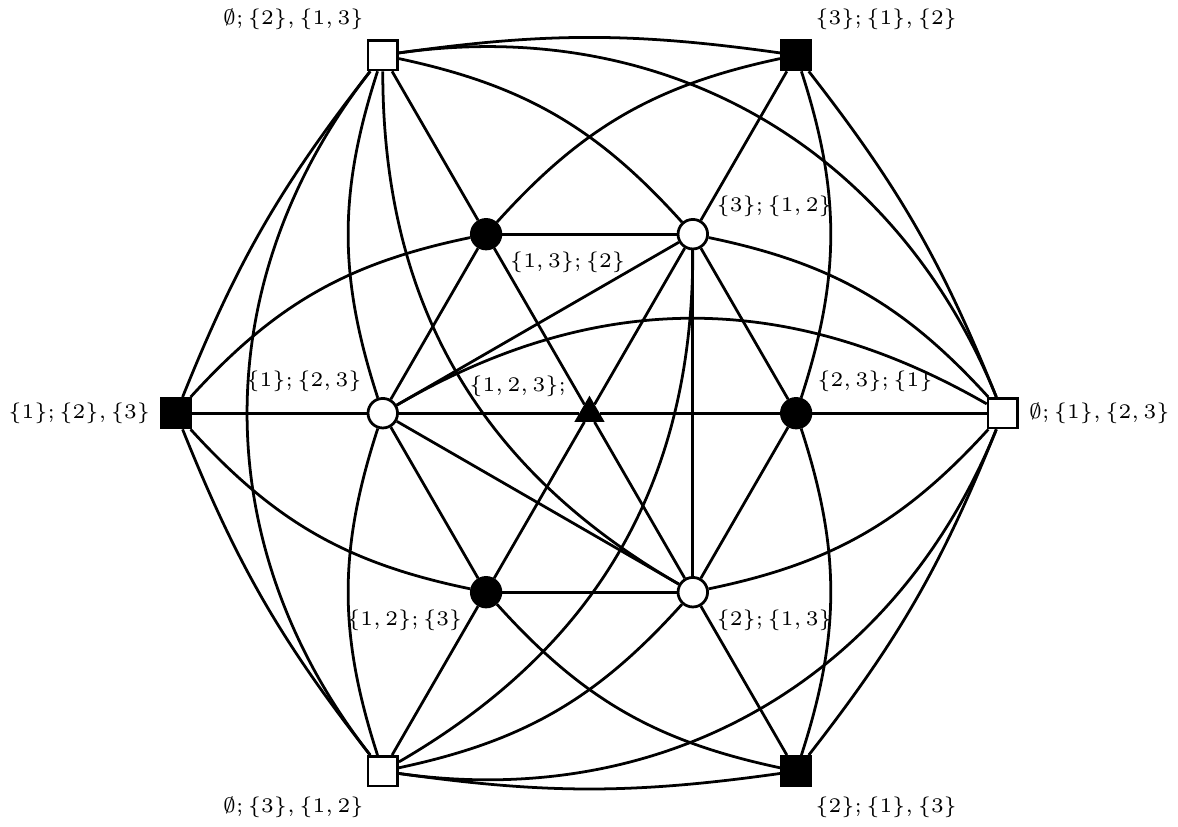}{0.7}{0.6}
        \end{minipage}
        \begin{minipage}{0.4\textwidth}
            \includefig{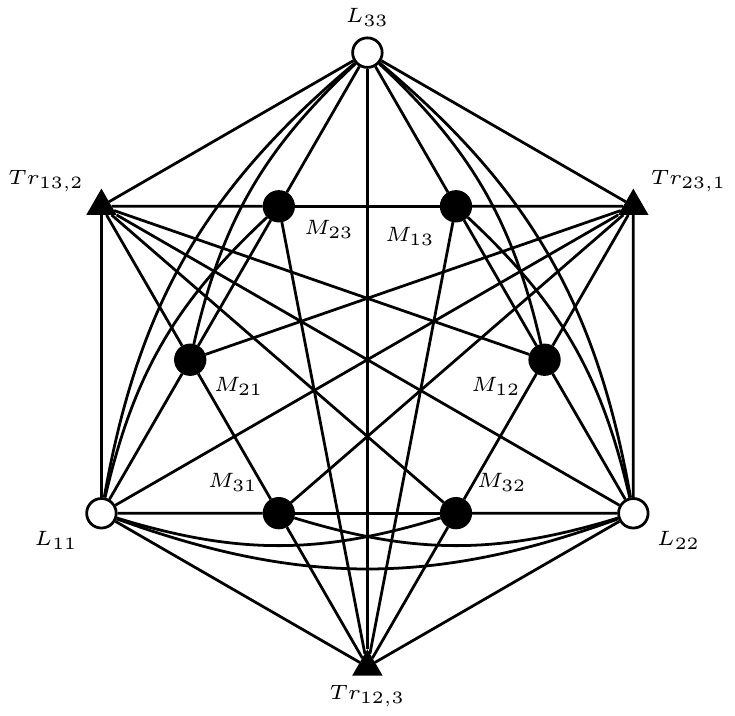}{0.6}{0.7}
        \end{minipage}
    }
    \caption{The skeleton and ridge graph of $\PMET_3 = 0,1$-$\PMET_3$}
    \label{fig:pmet3-ridge}
\end{figure}

\begin{conj}

(i) $\Sk(\OCUT_n)=K_{2S(n,2)}$ and belongs to
 $\Sk(\WQMET_n)$.

(ii) $\overline{\Sk(0,1\mbox{-}\dWMET_n)}=K_{1,S(n,2)}+S(n,2)K_2$;

 $\Sk(0,1$-$\dWMET_n)$ has diameter $2$,   all non-adjacencies are of the form:
 
  $(((0));(1))\nsim (\delta'(S);(0))$ and $(\delta'(S);w')\nsim
(\delta'(S);w')$.

\end{conj}

\begin{conj}

(i) $\Ri(\PMET_n)$ has diameter $2$,  all  non-adjacencies are:

$L_{ii}\nsim M_{ik}$; $M_{ij}\nsim M_{ji}, M_{ki}, M_{jk},Tr_{ij,k}$; $Tr_{ij,k}\nsim
Tr_{i'j',k'}$ if they  conflict.

(ii)  $\Ri(\WQMET_n)$ has diameter $2$; it is $\Ri(\PMET_n)$ without vertices $L_{ii}$.
\end{conj}

\clearpage

\begin{table}[t]\begin{center}
\scriptsize
\makebox[\textwidth]{
\begin{tabular}{|c|c|cccccc|c|c|c|}
 \hline
$R_i$ & Representative
 &1&21&2&31&32&3&Inc.&Adj.&$|R_i|$\\
 \hline
$R_1$ $\blacktriangle$&$(\delta(\emptyset);(1))$
&1&0&1&0&0&1&9&6&1\\
$R_2$ $\circ$&$(\delta(\{1\};w'')$ &0&1&1&1&0&1&9&8&3\\
$R_3$ $\bullet$&$(\delta(\{1\};w')$ &1&1&0&1&0&0&9&8&3\\
$R_4$ $\square$&$(\delta(\{1\};(0))$ &0&1&0&1&0&0&9&8&3\\
\hline \hline

$F_1$ $\circ$&$
L_1: w_1\ge 0$& 1&0&0&0&0&0&6&9&3\\
$F_2$ $\blacktriangle$&$Tr_{12,3}:  d_{13}+d_{23}-d_{12}\ge 0$&0&-1&0&1&1&0&7&8&3\\
$F_3$ $\bullet$&$
M'_{12}:  d_{12}+(w_2-w_1)
\ge 0$&-1&1&1&0&0&0&6&6&6\\
$F_4$ $\triangle$&$
Tr'_{12,3}: (d_{13}+d_{23}-d_{12})+2(w_1+w_2-w_3)\ge 
0$&2&-1&2&1&1&-2&5&5&3\\
 \hline \hline
\end{tabular}
}
\caption{The 
 orbits of extreme rays and facets in
$0,1$-$\dWMET_{3}$}
\label{tab:tabl0,1-3facets}
\end{center}\end{table}

\begin{figure}[t!]
	\vspace{-2mm}
    \centering
    \makebox[\textwidth]{
        \begin{minipage}{0.45\textwidth}
            \includefig{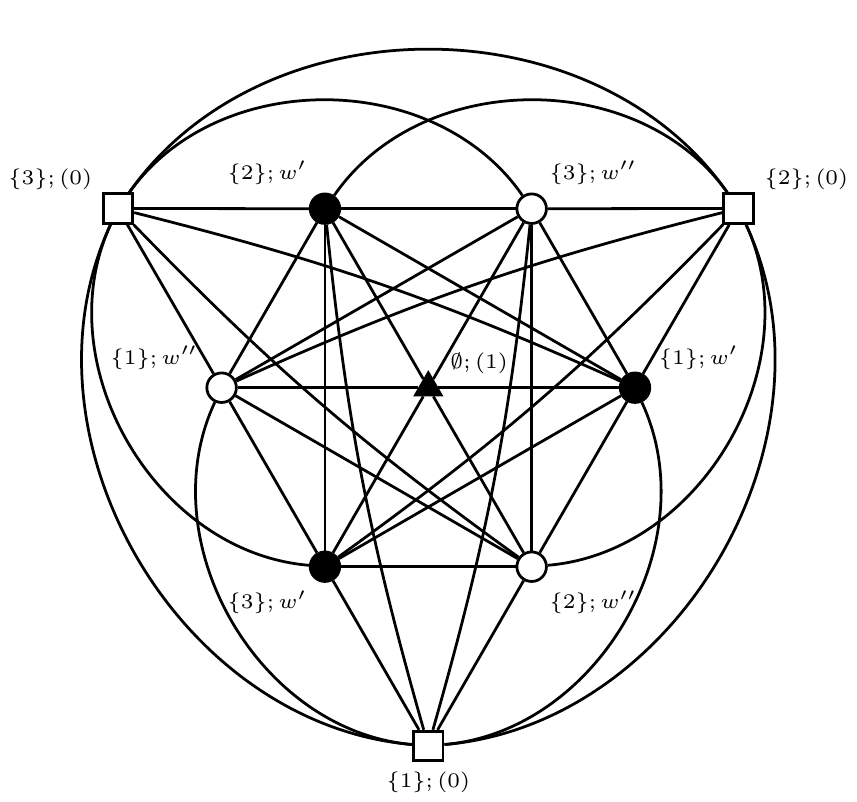}{0.7}{0.7}
        \end{minipage}
        \begin{minipage}{0.45\textwidth}
            \includefig{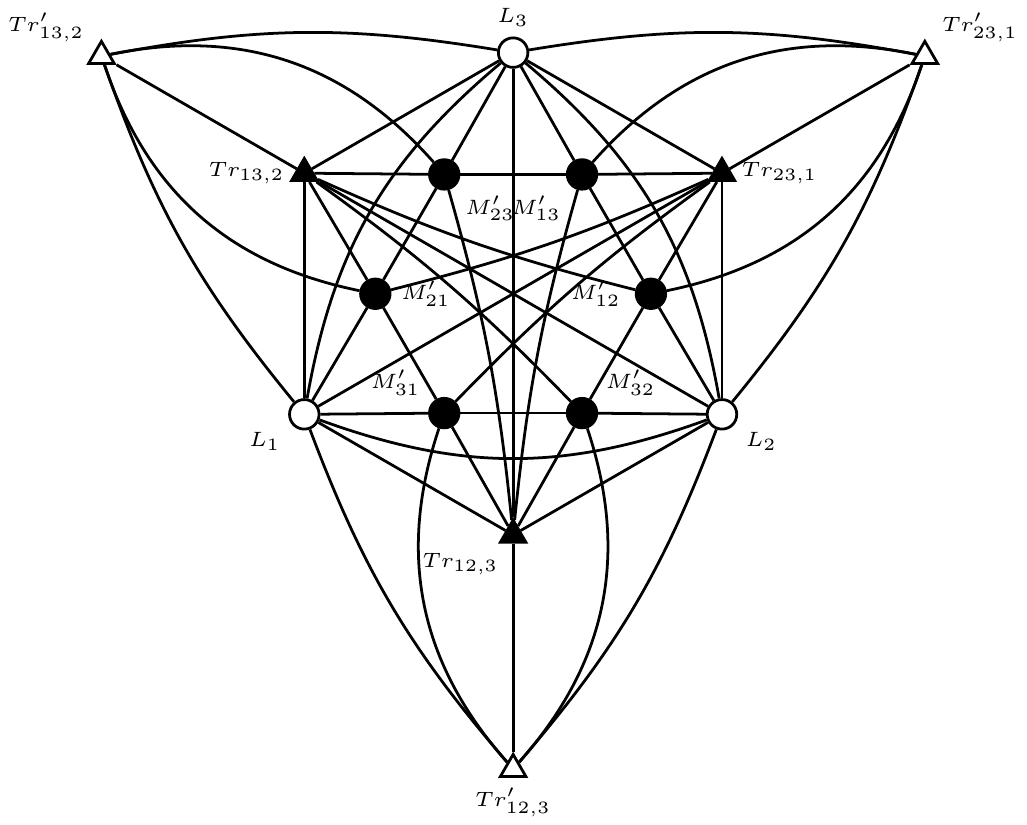}{0.6}{0.7}
        \end{minipage}
    }
    \caption{The skeleton and ridge graph of $0,1$-$\dWMET_3$}
    \label{fig:01dwmet3-skeleton}
\end{figure}

\end{document}